\documentclass[preprint,12pt]{elsarticle}

\usepackage{graphicx}

\usepackage[ruled,lined,algonl,boxed]{algorithm2e}

\usepackage{amssymb}
\usepackage{amsthm}
\newtheorem{thm}{Theorem}

\newtheorem{conj}[thm]{Conjecture}
\newtheorem{lemma}[thm]{Lemma}
\newdefinition{rmk}{Remark}
\newproof{pf}{Proof}
\newproof{pot}{Proof of Theorem \ref{thm2}}

\journal{\empty}

\begin{document}

\begin{frontmatter}



\title{Upper bounds for the bondage number of graphs on topological surfaces}


\author[label2]{Andrei Gagarin\corref{cor1}}
\ead{andrei.gagarin@acadiau.ca}
\cortext[cor1]{Corresponding author}
\author[label2a]{Vadim Zverovich}
\ead{vadim.zverovich@uwe.ac.uk}
\address[label2]{Department of Mathematics and Statistics, Acadia University, Wolfville, Nova Scotia, B4P 2R6, Canada}
\address[label2a]{Department of Mathematics and Statistics, University of the West of England, Bristol, BS16 1QY, UK}

\begin{abstract}
The bondage number $b(G)$ of a graph $G$ is the smallest number of edges of $G$ whose removal from $G$ results in a graph having the domination number larger than that of $G$. We show that, for a graph $G$ having the maximum vertex degree $\Delta(G)$ and embeddable on an orientable surface of genus $h$ and a non-orientable surface of genus $k$,
$$b(G)\le \min\{\Delta(G)+h+2,\ \Delta(G)+k+1\}.$$ This generalizes known upper bounds for planar and toroidal graphs.
\end{abstract}


\end{frontmatter}


\section{Introduction}
\label{Intro}
We consider simple finite non-empty graphs. 
For a graph $G$, its vertex and edge sets are denoted, respectively, by $V(G)$ and $E(G)$.
We also use the following standard notation: $d(v)$ for the degree of a vertex $v$ in $G$,
$\Delta=\Delta(G)$ for the maximum vertex degree of $G$, 
$\delta=\delta(G)$ for the minimum vertex degree of $G$, and 
$N(v)$ for the neighbourhood of a vertex $v$ in $G$. 

A set $D\subseteq V(G)$ is a {\it dominating set} if every vertex not in $D$ is
adjacent to at least one vertex in $D$. The minimum cardinality of
a dominating set of $G$ is the \emph{domination number}
$\gamma(G)$. Clearly, for any spanning subgraph $H$ of $G$, $\gamma(H)\ge\gamma(G)$. The \emph{bondage number} of $G$, denoted by $b(G)$, is the minimum cardinality of a set of edges $B\subseteq E(G)$ such that $\gamma(G-B)>\gamma(G)$.

The bondage number was introduced by Bauer et al. \cite{B83} (see also Fink et al. \cite{F90}). Two unsolved classical conjectures for the bondage number of arbitrary and planar graphs are as follows.
\begin{conj} [Teschner \cite{T95}] \label{all}
For any graph $G$, $b(G)\le \frac{3}{2}\Delta(G)$.
\end{conj}

Hartnell and Rall \cite{H94} and Teschner \cite{T96} showed that for the cartesian product $G_n = K_n \times K_n$, $n\ge2$, the bound of Conjecture \ref{all} is sharp, i.e. $b(G_n)= \frac{3}{2}\Delta(G_n)$. Teschner \cite{T95} also proved that Conjecture \ref{all} holds when $\gamma(G) \le 3$.

\begin{conj} [Dunbar et al. \cite{D98}] \label{planar-conj}
If $G$ is a planar graph, then $b(G)\le\Delta(G)+1$.
\end{conj}


The planar graphs are precisely the graphs that can be drawn on the sphere with no crossing edges.
A topological surface $S$ can be obtained from the sphere $S_0$ by adding a number of handles or crosscaps. If we add $h$ handles to $S_0$, we obtain an orientable surface $S_h$, which is often referred to as the \emph{$h$-holed torus}. The number $h$ is called the \emph{orientable genus} of $S_h$. If we add $k$ crosscaps to the sphere $S_0$, we obtain a non-orientable surface $N_k$. The number $k$ is called the \emph{non-orientable genus} of $N_k$. Any topological surface is homeomorphically equivalent either to $S_h$ ($h\ge 0$), or to $N_k$ ($k\ge 1$). For example, $S_1$, $N_1$, $N_2$ are the \emph{torus}, the \emph{projective plane}, and the \emph{Klein bottle}, respectively.

A graph $G$ is \emph{embeddable} on a topological surface $S$ if it admits a drawing on the surface with no crossing edges. Such a drawing of $G$ on the surface $S$ is called an \emph{embedding} of $G$ on $S$. Notice that there can be many different embeddings of the same graph $G$ on a particular surface $S$. The embeddings can be distinguished and classified by different properties. The set of faces of a particular embedding of $G$ on $S$ is denoted by $F(G)$.

An embedding of $G$ on the surface $S$ is a \emph{$2$-cell embedding} if each face of the embedding is homeomorphic to an open disk. In other words, a $2$-cell embedding is an embedding on $S$ that ``fits" the surface. This is expressed in Euler's formulae (\ref{orient}) and (\ref{non-orient}) of Theorem \ref{Euler}. For example, a cycle $C_n$ ($n\ge3$) does not have a $2$-cell embedding on the torus, but it has $2$-cell embeddings on the sphere and the projective plane. Similarly, a planar graph may have $2$-cell and non-$2$-cell embeddings on the torus.

The following result is usually known as (generalized) \emph{Euler's formula}. We state it here in a form similar to Thomassen \cite{T92}.
\begin{thm} [Euler's Formula, \cite{T92}] \label{Euler}
Suppose a connected graph $G$ with $|V(G)|$ vertices and $|E(G)|$ edges admits a $2$-cell embedding having $|F(G)|$ faces on a topological surface $S$. 
Then, either $S=S_h$ and
\begin{equation} \label{orient}
|V(G)| - |E(G)| + |F(G)| = 2-2h,
\end{equation}
or $S=N_k$ and
\begin{equation} \label{non-orient}
|V(G)| - |E(G)| + |F(G)| = 2-k.
\end{equation}
\end{thm}
Equation (\ref{orient}) is usually referred to as Euler's formula for an orientable surface $S_h$ of genus $h$, $h\ge 0$, and Equation (\ref{non-orient}) is known as Euler's formula for a non-orientable surface $N_k$ of genus $k$, $k\ge 1$.

The \emph{orientable genus} of a graph $G$ is the smallest integer $h=h(G)$ such that $G$ admits an embedding on an orientable topological surface $S$ of genus $h$.
The \emph{non-orientable genus} of $G$ is the smallest integer $k=k(G)$ such that $G$ can be embedded on a non-orientable topological surface $S$ of genus $k$. Clearly, in general, $h(G)\not=k(G)$, and the embeddings on $S_{h(G)}$ and $N_{k(G)}$ must be $2$-cell embeddings.

Trying to prove Conjecture \ref{planar-conj}, Kang and Yuan \cite{K00} came up with the following upper bound whose simpler topological proof was later discovered by Carlson and Develin \cite{C06}. 
\begin{thm} [\cite{K00,C06}] \label{planar}
For any connected planar graph $G$, 
$$b(G)\le \min\{8,\ \Delta(G)+2\}.$$
\end{thm}

This solves Conjecture \ref{planar-conj} in case $\Delta(G)\ge7$. The upper bound of Theorem \ref{planar} is for the sphere $S_0$ that has orientable genus $h=0$. 
The proof of Theorem \ref{planar} in \cite{C06} is topologically intuitive, uses Euler's formula for the sphere, and allows its authors to establish a partially similar result for the torus.
\begin{thm} [\cite{C06}] \label{toroidal}
For any connected toroidal graph $G$, $b(G)\le \Delta(G)+3$.
\end{thm}

Notice that the torus $S_1$ has orientable genus $h=1$. As mentioned in \cite{C06}, it is sufficient to prove the results of Theorems \ref{planar} and \ref{toroidal} for connected graphs because the bondage number of a disconnected graph $G$ is the minimum of the bondage numbers of its components.

In this paper, we prove the following result which generalizes the corresponding upper bounds of Theorems \ref{planar} and \ref{toroidal} for any orientable or non-orientable topological surface $S$.
\begin{thm} \label{any-surface}
For a connected graph $G$ of orientable genus $h$ and non-orientable genus $k$, 
$$b(G)\le \min\{\Delta(G)+h+2,\ \Delta(G)+k+1\}.$$
\end{thm}
The upper bound of Theorem \ref{any-surface} follows from Theorems \ref{thm-orient} and \ref{thm-non-orient} proved below in Section 2.


\section{The bondage number on orientable and non-orientable surfaces}
\label{Section2}
In this section, we prove Theorem \ref{any-surface} by considering orientable and non-orientable surfaces separately. The proofs are done by using Euler's formulae (\ref{orient}) and (\ref{non-orient}), counting arguments, and the following result.

\begin{lemma} [Hartnell and Rall \cite{H94}] \label{lem}
For any edge $uv$ in a graph $G$, we have
$b(G)\le d(u)+d(v)-1-|N(u)\cap N(v)|$. 
In particular, this implies that
$b(G) \le \delta(G)+\Delta(G)-1$ \emph{(see also \cite{B83,F90})}.
\end{lemma} 

Having a graph $G$ embedded on a surface $S$, each edge $e_i=uv \in E(G)$, $i=1,\ldots,|E(G)|$, can be assigned two weights, $w_i=\frac{1}{d(u)}+\frac{1}{d(v)}$ and $f_i=\frac{1}{m'}+\frac{1}{m''}$, where $m'$ is the number of edges on the boundary of a face on one side of $e_i$, and $m''$ is the number of edges on the boundary of the face on the other side of $e_i$. Notice that, in an embedding on a surface, an edge $e_i$ may be not separating two distinct faces, but instead it can appear twice on the boundary of the same face. For example, every edge of a path $P_n$ ($n\ge 2$) embedded on the sphere is on the boundary of a unique face, and it appears exactly twice on the face boundary walk: once for each side of the edge. Clearly, in this case, $m'=m''=2(n-1)$ and $f_i=\frac{2}{m'}=\frac{2}{m''}=\frac{1}{n-1}$. 

Notice that weights $w_i$ and $f_i$, $i=1,\ldots,|E(G)|$, count the number of vertices of $G$ and faces of its embedding on $S$ as follows: 
$$\sum_{i=1}^{|E(G)|}w_i=|V(G)|,\ \ \ \ \ \sum_{i=1}^{|E(G)|}f_i=|F(G)|.$$ 
Then, by Euler's formula (\ref{orient}), we have 
$$\sum_{i=1}^{|E(G)|}(w_i+f_i-1)=|V(G)|+|F(G)|-|E(G)|=2-2h,$$
or, in other words,
$$\sum_{i=1}^{|E(G)|}\left(w_i+f_i-1-\frac{2-2h}{|E(G)|}\right)=\sum_{i=1}^{|E(G)|}\left(w_i+f_i-1+\frac{2h-2}{|E(G)|}\right)=0.$$
Now, each edge $e_i=uv \in E(G)$, $i=1,\ldots,|E(G)|$, can be associated with the quantity $w_i+f_i-1+\frac{2h-2}{|E(G)|}$ called the \emph{oriented curvature} of the edge.
Also, by Euler's formula (\ref{non-orient}), we have 
$$\sum_{i=1}^{|E(G)|}(w_i+f_i-1)=|V(G)|+|F(G)|-|E(G)|=2-k,$$
or, in other words,
$$\sum_{i=1}^{|E(G)|}\left(w_i+f_i-1-\frac{2-k}{|E(G)|}\right)=\sum_{i=1}^{|E(G)|}\left(w_i+f_i-1+\frac{k-2}{|E(G)|}\right)=0.$$
Then, each edge $e_i=uv \in E(G)$, $i=1,\ldots,|E(G)|$, can be associated with the quantity $w_i+f_i-1+\frac{k-2}{|E(G)|}$ called the \emph{non-oriented curvature} of the edge.

\begin{thm} \label{thm-orient}
Let $G$ be a connected graph $2$-cell embeddable on an orientable surface of genus $h\ge 0$. Then
\begin{equation} \label{eqn-orient}
b(G)\le \Delta(G)+h+2.
\end{equation}
\end{thm}

\begin{pf} 
Suppose $G$ is $2$-cell embedded on the $h$-holed torus $S_h$. By Lemma \ref{lem}, if $G$ has any vertices of degree $h+3$ or less, we have $\delta(G)\le h+3$, and inequality (\ref{eqn-orient}) holds. Therefore, we can assume $\Delta(G)\ge \delta(G)\ge h+4$. 

Now, suppose the opposite, $b(G)\ge \Delta(G)+h+3$. Then, by Lemma \ref{lem}, for any edge $e_i=uv$, $i=1,\dots,|E(G)|$, we have
$d(u)+d(v)-1-|N(u)\cap N(v)|\ge b(G)\ge \Delta(G)+h+3$. Then, $d(u)+d(v)\ge \Delta(G)+h+4+|N(u)\cap N(v)|$, and $d(u)\le \Delta(G)$, $d(v)\le \Delta(G)$.
If either $d(u)$ or $d(v)$ is equal to $h+4$, the other degree must be equal to $\Delta(G)\ge h+4$, and $u$ and $v$ cannot have any common neighbors, so that $m'$ and $m''$ are both at least $4$. Since in this case $|E(G)|\ge \frac{(h+4)(h+5)}{2}$, such an edge $e_i=uv$ has a negative oriented curvature:
$$w_i+f_i-1+\frac{2h-2}{|E(G)|}\le \frac{2}{h+4}+\frac{2}{4}-1+\frac{2(2h-2)}{(h+4)(h+5)}=\frac{-8+h(3-h)}{2(h+4)(h+5)}<0$$
for any $h\ge 1$, and, in case $h=0$, 
$$w_i+f_i-1-\frac{2}{|E(G)|}\le \frac{1}{4}+\frac{1}{4}+\frac{1}{4}+\frac{1}{4}-1-\frac{2}{|E(G)|}=\frac{-2}{|E(G)|}<0.$$

Suppose one of $d(u)$ and $d(v)$ is equal to $h+5$, without loss of generality, $d(u)=h+5$. Then, 
$\Delta(G)\ge d(v)\ge \Delta(G)-1+|N(u)\cap N(v)|$.
If $d(v)=h+4=\Delta(G)-1$, we are in the previous case. Otherwise, we have $d(v)\ge h+5$, and at most one of $m'$ and $m''$ can be equal to $3$, implying the other is at least $4$. Then again, since in this case $|E(G)|\ge \frac{(h+4)(h+4)+2(h+5)}{2}=\frac{h^2+10h+26}{2}$, the edge $e_i$ must have a negative oriented curvature:
$$w_i+f_i-1+\frac{2h-2}{|E(G)|}\le \frac{2}{h+5}+\frac{1}{3}+\frac{1}{4}-1+\frac{2(2h-2)}{h^2+10h+26}=\frac{-5h^3-3h^2+52h-266}{12(h+5)(h^2+10h+26)}<0$$
for any $h\ge 1$, and, in case $h=0$,
$$w_i+f_i-1-\frac{2}{|E(G)|}\le \frac{1}{5}+\frac{1}{5}+\frac{1}{3}+\frac{1}{4}-1-\frac{2}{|E(G)|}=-\frac{1}{60}-\frac{2}{|E(G)|}<0.$$

The only remaining case is when $d(u)\ge h+6$ and $d(v)\ge h+6$. Since $m'\ge 3$ and $m''\ge 3$, and, in this case, $|E(G)|\ge \frac{(h+4)(h+5)+2(h+6)}{2}=\frac{h^2+11h+32}{2}$, the edge $e_i$ must have a negative oriented curvature:
$$w_i+f_i-1+\frac{2h-2}{|E(G)|}\le \frac{2}{h+6}+\frac{2}{3}-1+\frac{2(2h-2)}{h^2+11h+32}=\frac{-h^3+h^2+28h-72}{3(h+6)(h^2+11h+32)}<0$$
for any $h\ge 1$, and, in case $h=0$,
$$w_i+f_i-1-\frac{2}{|E(G)|}\le \frac{1}{6}+\frac{1}{6}+\frac{1}{3}+\frac{1}{3}-1-\frac{2}{|E(G)|}=\frac{-2}{|E(G)|}<0.$$

Summing over all edges $e_i\in E(G)$ yields
$$\sum_{i=1}^{|E(G)|}\left(w_i+f_i-1+\frac{2h-2}{|E(G)|}\right)<0,$$
which is a contradiction to Euler's formula (\ref{orient}) stating
$$\sum_{i=1}^{|E(G)|}\left(w_i+f_i-1-\frac{2-2h}{|E(G)|}\right)=|V(G)|+|F(G)|-|E(G)|-(2-2h)=0.$$
Thus, $b(G)\le \Delta(G)+h+2$.
\qed
\end{pf}


\begin{thm} \label{thm-non-orient}
Let $G$ be a connected graph $2$-cell embeddable on a non-orientable surface of genus $k\ge 1$. Then
\begin{equation} \label{eqn-non-orient}
b(G)\le \Delta(G)+k+1.
\end{equation}
\end{thm}

\begin{pf}
Suppose $G$ is $2$-cell embedded on the sphere with $k$ crosscaps $N_k$. By Lemma \ref{lem}, if $G$ has any vertices of degree $k+2$ or less, we have $\delta(G)\le k+2$, and inequality (\ref{eqn-non-orient}) holds. Therefore, we can assume $\Delta(G)\ge \delta(G)\ge k+3$.

Suppose the opposite, $b(G)\ge \Delta(G)+k+2$. 
Then, by Lemma \ref{lem},  for any edge $e_i = uv$, $i = 1,\ldots,|E(G)|$, we have 
$d(u)+d(v)-1-|N(u)\cap N(v)|\ge b(G)\ge \Delta(G)+k+2$. Then, $d(u)+d(v)\ge \Delta(G)+k+3+|N(u)\cap N(v)|$, and $d(u)\le \Delta(G)$, $d(v)\le \Delta(G)$.
If either $d(u)$ or $d(v)$ is equal to $k+3$, the other degree must be equal to $\Delta(G)\ge k+3$, and $u$ and $v$ cannot have any common neighbors, so that $m'$ and $m''$ are both at least $4$. Since in this case $|E(G)|\ge \frac{(k+3)(k+4)}{2}$, the non-oriented curvature of the edge $e_i=uv$ is
$$w_i+f_i-1+\frac{k-2}{|E(G)|}\le \frac{2}{k+3}+\frac{2}{4}-1+\frac{2(k-2)}{(k+3)(k+4)}=\frac{-4+k(1-k)}{2(k+3)(k+4)}<0$$
for any $k\ge 2$, and, in case $k=1$,
$$w_i+f_i-1-\frac{1}{|E(G)|}\le \frac{1}{4}+\frac{1}{4}+\frac{1}{4}+\frac{1}{4}-1-\frac{1}{|E(G)|}=\frac{-1}{|E(G)|}<0.$$

Suppose one of $d(u)$ and $d(v)$, let us say $d(u)$, is equal to $k+4$. 
Then, $\Delta(G)\ge d(v)\ge \Delta(G)-1+|N(u)\cap N(v)|$.
If $d(v)=k+3=\Delta(G)-1$, we are in the previous case. Otherwise, we have $d(v)\ge k+4$, and at most one of $m'$ and $m''$ can be equal to $3$, implying the other is at least $4$. Then again, since in this case $|E(G)|\ge\frac{(k+3)(k+3)+2(k+4)}{2}=\frac{k^2+8k+17}{2}$, the edge $e_i$ must have a negative non-oriented curvature:
$$w_i+f_i-1+\frac{k-2}{|E(G)|}\le \frac{2}{k+4}+\frac{1}{3}+\frac{1}{4}-1+\frac{2(k-2)}{k^2+8k+17}=\frac{-124-5k-12k^2-5k^3}{12(k+4)(k^2+8k+17)}<0$$
for any $k\ge 2$, and, in case $k=1$,
$$w_i+f_i-1-\frac{1}{|E(G)|}\le \frac{1}{5}+\frac{1}{5}+\frac{1}{3}+\frac{1}{4}-1-\frac{1}{|E(G)|}=-\frac{1}{60}-\frac{1}{|E(G)|}<0.$$

The only remaining case is when $d(u)\ge k+5$ and $d(v)\ge k+5$. Since $m'\ge 3$ and $m''\ge 3$, and, in this case, $|E(G)|\ge \frac{(k+3)(k+4)+2(k+5)}{2}=\frac{k^2+9k+22}{2}$, the edge $e_i$ must have a negative non-oriented curvature:
$$w_i+f_i-1+\frac{k-2}{|E(G)|}\le \frac{2}{k+5}+\frac{2}{3}-1+\frac{2(k-2)}{k^2+9k+22}=\frac{-k^3-2k^2+5k-38}{3(k+5)(k^2+9k+22)}<0$$
for any $k\ge 2$, and, in case $k=1$,
$$w_i+f_i-1-\frac{1}{|E(G)|}\le \frac{1}{6}+\frac{1}{6}+\frac{1}{3}+\frac{1}{3}-1-\frac{1}{|E(G)|}=\frac{-1}{|E(G)|}<0.$$

Summing over all edges $e_i\in E(G)$ yields
$$\sum_{i=1}^{|E(G)|}\left(w_i+f_i-1+\frac{k-2}{|E(G)|}\right)<0,$$
which is a contradiction to Euler's formula (\ref{non-orient}) stating
$$\sum_{i=1}^{|E(G)|}\left(w_i+f_i-1-\frac{2-k}{|E(G)|}\right)=|V(G)|+|F(G)|-|E(G)|-(2-k)=0.$$
Thus, $b(G)\le \Delta(G)+k+1$, and the proof is complete.
\qed
\end{pf}

\section{Conclusions}
\label{FR}
The upper bound of Theorem \ref{any-surface} provides a hierarchy of upper bounds that eventually may help solving Conjecture \ref{all}. However, it can be seen that the bounds of Theorems \ref{thm-orient} and \ref{thm-non-orient} are not tight for larger values of the genera $h=h(G)$ and $k=k(G)$. For example, by adjusting respectively the proofs of Theorems \ref{thm-orient} and \ref{thm-non-orient}, upper bound (\ref{eqn-orient}) can be improved to $b(G)\le \Delta(G)+h+1$ for $h\ge 8$, and upper bound (\ref{eqn-non-orient}) can be improved to $b(G)\le \Delta(G)+k$ for $k\ge 3$ and to $b(G)\le \Delta(G)+k-1$ for $k\ge 6$. It is left to the reader to adjust the proofs and bounds for a particular topological surface of higher genus. 

The bounds of Theorems \ref{thm-orient} and \ref{thm-non-orient} are presented in this form for clarity and simplicity of presentation and proofs. In general, one may try to find certain (linear or sublinear) functions of $h$ and $k$ to improve the bounds by replacing the terms $h+2$ and $k+1$, respectively, or to provide asymptotically better bounds. For example, simple asymptotic improvements follow from the upper bounds on the minimum vertex degree of graphs embeddable on topological surfaces: $\delta(G)\le\lfloor\frac{5+\sqrt{1+48h}}{2}\rfloor$ for $h\ge 1$, 
$\delta(G)\le\lfloor\frac{5+\sqrt{1+24k}}{2}\rfloor$ for $k\ge 2$ (e.g., see Sachs \cite{Sachs}), and, for a planar or projective-planar graph $G$, i.e. when $h=0$ or $k=1$, $\delta(G)\le 5$. Then, from Lemma \ref{lem}, $b(G)\le \Delta(G)+\lfloor\frac{3+\sqrt{1+48h}}{2}\rfloor$ for $h\ge1$ and $b(G)\le \Delta(G)+\lfloor\frac{3+\sqrt{1+24k}}{2}\rfloor$ for $k\ge1$, which are better than bounds (\ref{eqn-orient}) and (\ref{eqn-non-orient}), respectively, when $h\ge 12$ and $k\ge 8$. However, for example, an adjusted proof of Theorem \ref{thm-non-orient} gives $b(G)\le \Delta +k-411=\Delta+53$ for $k=464$, which is better than $b(G)\le \Delta(G)+\lfloor\frac{3+\sqrt{1+24k}}{2}\rfloor=\Delta +54$ in this case. Therefore, adjustments of the proofs of Theorems \ref{thm-orient} and \ref{thm-non-orient} can provide better results than some asymptotic improvements. However, it would be interesting to have an asymptotic improvement that provides a certain justification of its quality.


In view of Theorem \ref{planar}, its proof in \cite{C06}, and results presented in this paper, it should be reasonable to conjecture that, when $\Delta(G)$ is sufficiently large, the bondage number $b(G)$ is bounded by a certain constant depending only on the properties of topological surfaces where $G$ embeds.
\begin{conj}
For a connected graph $G$ of orientable genus $h$ and non-orientable genus $k$, $b(G)\le\min\{c_h,\, c^\prime_k,\, \Delta(G)+h+2,\, \Delta(G)+k+1\}$, where $c_h$ and $c^\prime_k$ are constants depending, respectively, on the orientable and non-orientable genera of $G$.
\end{conj}

Since $\delta(G)\le5$ for a planar graph $G$, Fischermann et al. \cite{F03} ask whether there exist planar graphs of bondage numbers $6$, $7$, or $8$. A class of planar graphs with the bondage number equal to $6$ is shown in \cite{C06}. Therefore, in case of planar graphs, we have $6\le c_0\le 8$. It would be interesting to have an estimation for the constants $c_h$ and $c^\prime_k$ for the torus $S_1$, projective plane $N_1$, and Klein bottle $N_2$.





\end{document}